\documentclass[twocolumn]{jsiamletters}
\group{Topological Data Analysis}
\affiliation{
  Center for Artificial Intelligence and Mathematical Data Science (Angels), Okayama University
}{
  3-1-1 Tsushima-naka, Kita-ku, Okayama 700-8530 Japan
}
\affiliation{
  Photon Factory, Institute of Materials Structure Science, High Energy Accelerator Research Organization (KEK)
}{
  Tsukuba, Ibaraki 305-0801, Japan
}
\affiliation{
Department of Materials Structure Science, School of High Energy Accelerator Science, SOKENDAI (The Graduate University for Advanced Studies)
}{
  Tsukuba, Ibaraki 305-0801, Japan
}

\authorinfo{Ippei Obayashi}{1*}{i.obayashi@okayama-u.ac.jp}
\authorinfo{Masao Kimura}{2,3}{masao.kimura@kek.jp}
\email{i.obayashi@okayama-u.ac.jp}
\title{Persistent homology analysis with nonnegative matrix factorization for 3D voxel data of iron ore sinters}
\abstract{
  This paper proposes a data analysis method using persistent homology and nonnegative matrix factorization.
  A concatenated persistence image technique is used to extract
  coexisting structures from the persistence diagrams of different dimensions hidden behind the data.
  To demonstrate the potential of our method, we apply the method to 3D voxel data of iron ore sinters obtained by X-ray computed tomography.
  The analysis successfully captures the coexistence structures in these iron ore sinters.  
}
\newcommand{\pd}{D}
\newcommand{\R}{\mathbb{R}}
\newcommand{\cfd}{\mathrm{CFD}}

\begin{document}

\maketitle

\section{Introduction}
Persistent homology (PH) is a data analysis framework that uses topology, which enables us to extract quantitative information about the shape of the data~\cite{homcloud}. PH is good at handling disordered and heterogeneous structures and has been applied to various fields, including materials science and life science.

In PH, the data are converted into a two-dimensional (2D) histogram called a \emph{persistence diagram} (PD) and then analysis is carried out using this diagram. If the PD has some outstanding features, the analysis is easy. In reality, however, we often encounter cases where such features cannot easily be found. One way to deal with this problem is to apply machine learning methods to the PDs. With this method, information about the shape of the data is quantified with PH, and then machine learning is used to find any characteristic patterns.

There are already many studies on the combination of PH and machine learning\cite{phml}.
There have also been various studies on the vectorization of PDs and kernel methods to use them as inputs for machine learning methods.
In this letter, we will present a new method that uses a vectorization method called persistence images (PIs)\cite{pi} and a machine learning method called nonnegative matrix factorization (NMF)\cite{nmf}. 
By concatenating PI vectors before applying NMF, we can successfully capture the coexistence structures shown in different dimensional PDs.
The concatenation technique is simple but works well in combination with NMF.

We also show a case study of analyzing three-dimensional (3D) X-ray computed tomography (X-CT) images of iron ore sinters using the proposed method. Iron ore sinters are the starting material for iron-making processes and is produced by agglomerating fine-size materials, including iron ore particles, limestone flux, and coke breeze through a sintering process. Then, the sinters are reduced in blast furnaces into pig iron \cite{iron}. The reduction reaction progresses heterogeneously due to the existence of multiple phases and pore networks in the sinters, and can result in crack formation, one of the main factors in sinter degradation. Characterization of the changes in the distribution of phases during the reduction reactions by the proposed method is expected to provide indispensable information for controlling the microstructures of iron ore sinters.

\section{Method}

\begin{figure*}[bt]
  \centering
  \includegraphics[width=0.8\hsize]{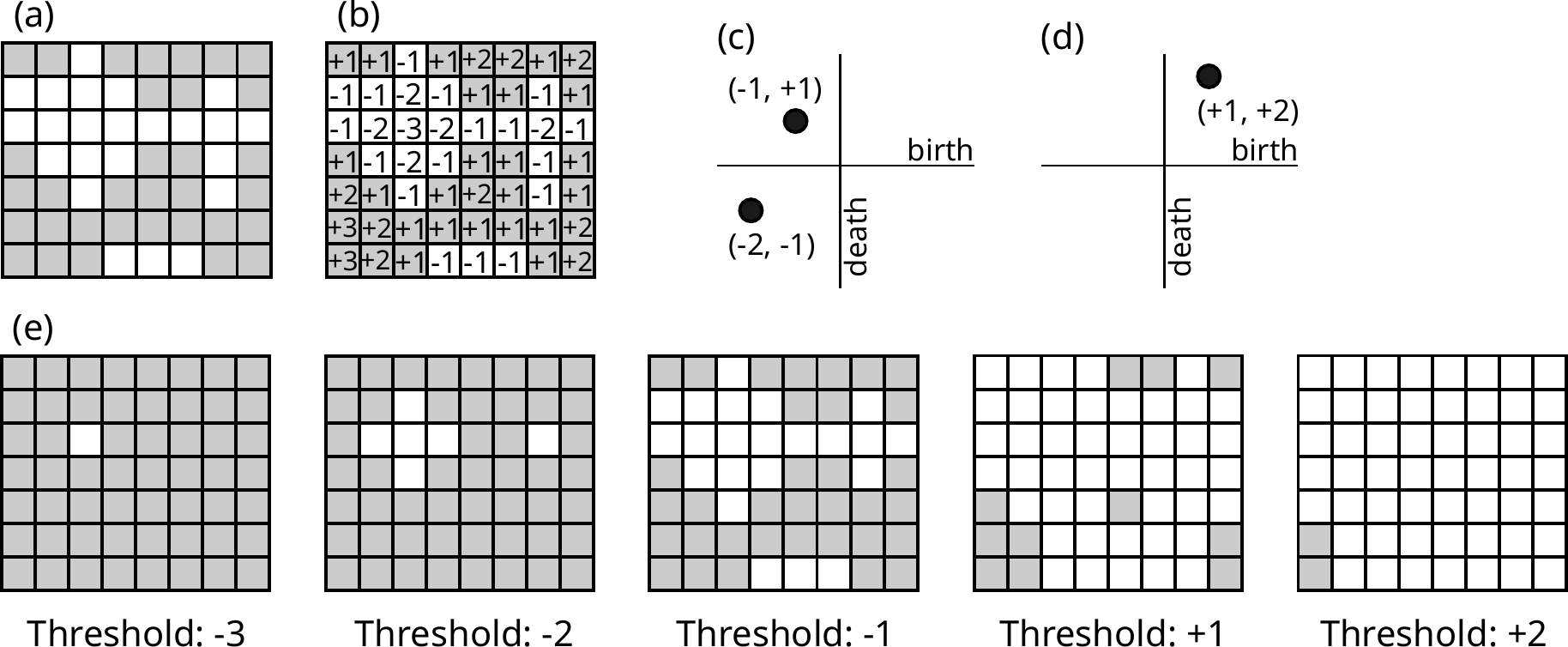}
  \caption{Persistent homology. (a) Input binary image. (b) Signed distance transform with Manhattan distance. (c) 0-dimensional PD. (d) 1-dimemsional PD.
    (f) Filtration of a binary image with respect to the signed distance transform.}
  \label{fig:ph}
\end{figure*}

To capture the geometric features of image data, we consider the change of homology in an increasing sequence called a \emph{filtration}.
We often use the \emph{signed distance transform} to build a filtration from a binary image.
Figure~\ref{fig:ph}(b) is the signed distance transform of Fig. \ref{fig:ph}(a).
For each gray pixel, its pixel value is the distance to the nearest white pixel. For each white pixel, its pixel value is the negative distance to the nearest gray pixel. In this paper, we use the Manhattan distance.
The filtration Fig.~\ref{fig:ph}(e) is obtained by changing the binarization threshold in Fig.~\ref{fig:ph}(b).
In this filtration, homology generators appear and disappear.
The theory of PH ensures a unique pairing of the appearance and disappearance of homology generators in the filtration.
The pair is called a \emph{birth-death pair}, and the collection of birth-death pairs is called a PD.
A PD is often visualized with a 2D scatter plot or a histogram, as shown in Fig.~\ref{fig:ph}(c)(d).

For the given voxel data $\{X_n\}_n$, $k$-dimensional PDs $\{\pd_k(X_n)\}$ for $k=0, 1, 2$ are computed using the signed distance transform.
From a PD $\pd$, we define the \emph{persistence surface} $\rho_D$ on $\R^2$ as follows:
\begin{equation}
  \begin{aligned}
    \rho_D(x, y) &= \sum_{(b, d) \in \pd} w(b, d) \exp \frac{(x-b)^2 + (y-d)^2}{2\sigma^2}, \\
    w(b, d) &= \arctan C(b-d)^p,
  \end{aligned}
\end{equation}
where $\sigma > 0$, $C > 0$, and $p \geq 1$ are parameters.
The distribution $\rho_D$ is discretized by a finite grid $G$, and we can regard the result as a vector.
The dimension of this vector is the same as the number of bins in the grid $G$.
The discretized vector $\hat{\rho}_D$ is called a \emph{persistence image} (PI).
Intuitively, the PI is the histogram of a PD over grid $G$.
Gaussian smoothing is applied to the histogram, and $w(b, d)$ is used to reflect the importance of points in $D$ (a birth-death pair near the diagonal is known to be less important than a pair far from the diagonal).
One of the advantages of PIs is that vectors of the same dimension can be converted to one PD histogram.

We apply \emph{nonnegative matrix factorization} (NMF) to the PIs.
From nonnegative vectors $\{v_1, \ldots, v_N\}$, NMF finds a small number of nonnegative vectors $\{w_1, \ldots w_M\}$ and nonnegative numbers $\{\lambda_{n,m}\}_{1 \leq n \leq N, 1 \leq m \leq M}$ that satisfy the following approximation in $\ell^2$-norm:
\begin{equation}
  v_n \approx \sum_{m=1}^M \lambda_{n, m}w_m.
\end{equation}
By applying NMF to PIs, we obtain the following approximation:
\begin{equation}
  \rho_{D_k(X_n)} \approx \sum_{m=1}^M \lambda_{n,m} \rho_m,
\end{equation}
where $\{\lambda_{m,n}\}$ are nonnegative numbers and $\{\rho_m\}_{m=1}^M$ are distributions on the $xy$-plane.
In this paper, we call $\lambda_{n,m}$ an \emph{NMF coefficient} and $\rho_m$ a \emph{feature distribution}.
Here, we utilize the fact that PIs are positive vectors.

To integrate PDs of different dimensions, we propose applying NMF to the concatenated vectors.
By concatenating three PIs $\hat{\rho}_{D_0(X_n)}$, $\hat{\rho}_{D_1(X_n)}$, and $\hat{\rho}_{D_2(X_n)}$ vertically and applying NMF to the concatenated vectors, we obtain the following approximation:
\begin{equation}\label{eq:approx-nmf3}
  \begin{bmatrix}
    \rho_{D_0(X_n)} \\ \rho_{D_1(X_n)} \\ \rho_{D_2(X_n)}
  \end{bmatrix}
  \approx
  \sum_{m=1}^M \lambda_{n,m}
  \begin{bmatrix}
    \rho_{0, m} \\ \rho_{1, m} \\ \rho_{2, m}
  \end{bmatrix}.
\end{equation}
In this approximation, we can regard the triple $[\rho_{0,m}, \rho_{1,m}, \rho_{2,m}]$ as a coexisting structure in the PDs of different dimensions.
We call $[\rho_{0,m}, \rho_{1,m}, \rho_{2,m}]$ the $m$th \emph{concatenated feature distribution} and refer to it as $\cfd_m$.

By analyzing $\rho_{d,m}$ and $\lambda_{n,m}$, we can understand which birth-death pairs are important for characterizing the PDs.
For a deeper investigation of such birth-death pairs, we apply ``inverse analysis'' to identify the geometric origins of the birth-death pairs\cite{inverse}.

The advantages of our method are summarized as follows:
\begin{itemize}
\item NMF is a technique similar to principal component analysis (PCA), and the advantages of NMF compared with PCA are as follows.
  These are advantageous for the analysis of PDs:
  \begin{itemize}
  \item Better interpretability because of the positivity of $\rho_m$ and $\lambda_{m,n}$.
  \item Each component of NMF tends to reflect the cluster structure of the data\cite{nmf-cluster}. This fact also contributes to the interpretability of the result.
  \end{itemize}
\item The technique using concatenated PIs extracts the coexisting structures in the PDs of different dimensions.
\item The fundamental idea is simple. Therefore, we can understand and implement the idea easily.
\end{itemize}
The presented analysis framework is similar to our previous work\cite{ph-linear}, but the introduction of NMF with concatenated PIs improves the interpretability of the result.

\section{Result and discussion}

We applied the proposed method to X-CT 3D voxel images of iron ore sinters.
We used the same data in previous research\cite{iron}.
The readers can see detailed information about the data (experimental setup, measurement devices, resolutions of the voxel images, etc.) in that paper.

In this experiment, six specimens were measured: two from the early stage of the reduction reaction, two more from the intermediate stage, and the last two from the final stage.
The specimens have mainly three phases: iron oxides, calcium ferrites, and large voids.
In the preprocessing step, the voxel images corresponding to iron oxides were extracted using an image segmentation technique (the watershed algorithm was used) after applying a Gaussian denoising filter.
After binarization, the voxel data were divided into $150 \times 150 \times 150$ cubes to investigate the relative local structures. The cubes corresponding to large initial pores (more than 40\% in the cube) were removed since such cubes were considered not to have useful information.
PDs were computed from the cubes using HomCloud\cite{homcloud}(\url{https://homcloud.dev/}). 
PIs were also computed using HomCloud, and NMF was performed using scikit-learn\cite{sklearn}.
The number of NMF components, $M$, was set to 3.
We tried different vectorization parameters and NMF parameters to check the robustness of the analysis, and we confirmed that the results do not change much even if there are small parameter changes.

\begin{figure}
  \centering
  \includegraphics[width=0.9\hsize]{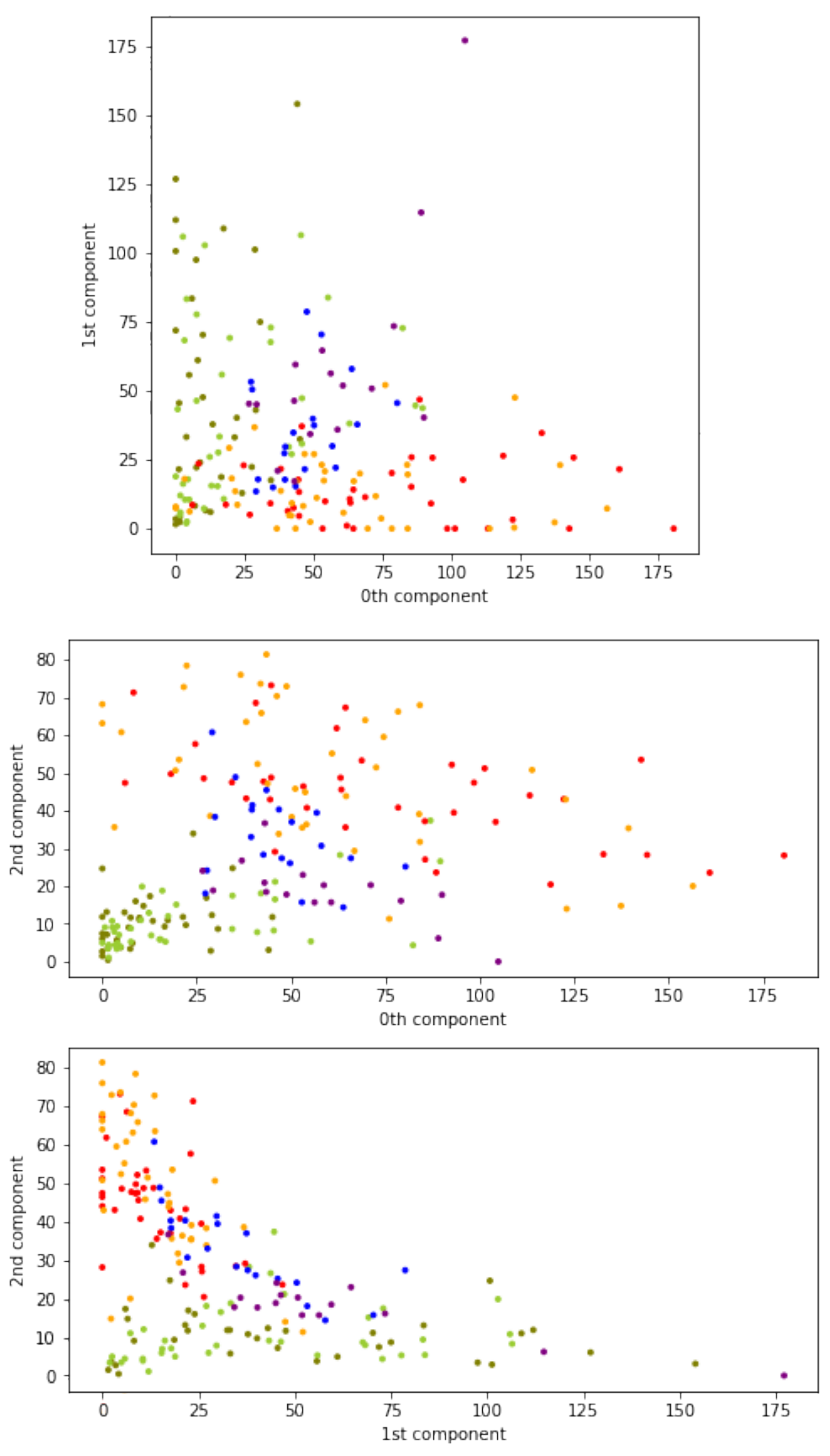}
  \caption{NMF coefficients: the color indicates the stage of reduction (green and olive for the early stage, red and orange for the intermediate stage, and blue and purple for the final stage). (a) $(\lambda_{*, 0}, \lambda_{*, 1})$ plot. (b) $(\lambda_{*, 0}, \lambda_{*, 2})$ plot. (c) $(\lambda_{*, 1}, \lambda_{*, 2})$ plot.}
  \label{fig:nmf-proj}
\end{figure}


Figure~\ref{fig:nmf-proj} shows the 3D NMF projection $\{(\lambda_{n, 0},\lambda_{n, 1},\lambda_{n, 2})\}$ as 2D projections.
Each point on the figure corresponds to one cube in the data, and the color indicates the stage of reduction (green and olive for the early stage, red and orange for the intermediate stage, and blue and purple for the final stage).

First, Fig.~\ref{fig:nmf-proj} shows that the features extracted by NMF are similar between different specimens at the same reaction stage.
This means that our method successfully detects differences between reaction stages.
Furthermore, the results indicate that the geometrical features of intermediate stage data mainly appear in the zeroth and second concatenated feature distributions and that the geometrical features of the early stage mainly appear in the first concatenated feature distribution.
The results also show that the data from the final stage have features between the early and intermediate stage data.


\begin{figure}
  \centering 
  \includegraphics[width=\hsize]{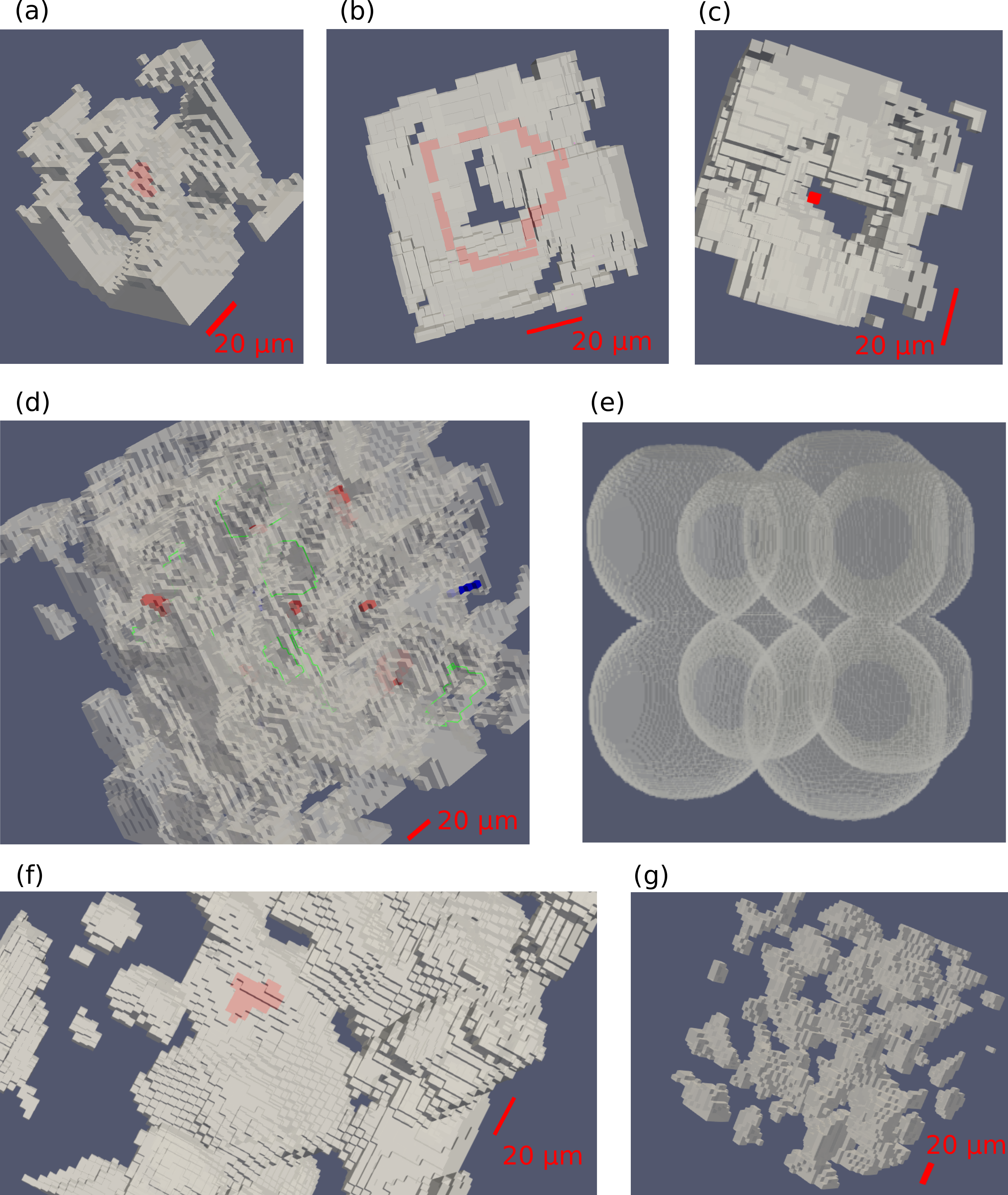}
  \caption{Geometric origins of feature distributions. (a)(b)(c) $\rho_{0,0}$, $\rho_{1,0}$, and $\rho_{2, 0}$, respectively, in $\cfd_0$.
    (d) Three types of structures corresponding to $\rho_{0,0}$ (red), $\rho_{1,0}$ (green), and $\rho_{2, 0}$ (blue).
    (e) Schematic picture of the coexisting structures corresponding to $\cfd_0$.
    (f) $\cfd_1$ (red voxels indicate the center of the structure).
    (g) $\cfd_2$.
  }
  \label{fig:ia}
\end{figure}

By applying inverse analysis to the result, we can identify the geometric origins of $\rho_{i,j}$.
Figure~\ref{fig:ia} shows some typical examples of the results. $\rho_{0, 0}$ corresponds to branched structures (Fig.~\ref{fig:ia}(a)), $\rho_{1, 0}$ corresponds to tunnels (Fig.~\ref{fig:ia}(b)), and $\rho_{2, 0}$ corresponds to hollow areas (Fig.~\ref{fig:ia}(c)).
The NMF results suggest that these three types of structures tend to coexist.
Figure~\ref{fig:ia}(d) shows an example of the three types of structures coexisting.

Figure~\ref{fig:ia}(e) shows the schematic picture of the coexisting structures corresponding to $\cfd_0$ as suggested by the obtained results.
In the figure, eight ball shapes can be seen very closely together.
$\rho_{0, 0}$ corresponds to the center of each ball, $\rho_{1, 0}$ corresponds to the tunnels in the figure, and $\rho_{2, 0}$ corresponds to the cavities surrounded by the balls.a
From the discussion, we conclude that the tunnel structure of iron oxide with complex branching is considered to correspond to the zeroth concatenated feature distribution.
In other words, one typical microstructure of iron ore sinters can be successfully described in terms of the zeroth concatenated feature distribution.
We note that the tunnels are actually filled with calcium ferrites and are not void.

In the same way, the geometric origins corresponding to $\cfd_1$ and $\cfd_2$ can be obtained. It can be considered that $\cfd_1$ (typical of the early stage) corresponds to relatively large grains of iron oxide with bumps (Fig.~\ref{fig:ia}(f)) and that $\cfd_2$ (another structure, typical of the intermediate stage) corresponds to the cluster of small grains of iron oxide (Fig.~\ref{fig:ia}(g)).

The proposed method of the analysis showed that complicated 3D microstructures found in iron ore sinters during the reduction reactions, which previously were not possible to understand quantitatively, can be successfully categorized and quantified. 
The three typical structures described above can be determined to be the most important shapes that appear during the reduction reactions even without any knowledge of the materials science of the reactions.
It should be noted that the determined structures are indeed consistent with our previous researchin terms of the materials science \cite{iron}. 
The structures corresponding to $\cfd_1$ dominate in number in the early stage, and as the reaction proceeds, their numbers decrease, while structures corresponding to $\cfd_0$ and $\cfd_2$ increase in number. At the final stage, structures corresponding to $\cfd_0$ and $\cfd_2$ decrease in number.

\section{Concluding remarks}

We now consider a scenario consistent with the results from the viewpoint of materials science.
The relatively large grains corresponding to $\cfd_1$ are likely to be iron ore particles.
In the reduction reaction process, various reactions occur heterogeneously due to multiple phases and the presence of a pore network, and the tunnel structures corresponding to $\cfd_0$ are thus formed from the reaction of relatively large grains, while the small grains corresponding to $\cfd_2$ are formed by the decomposition of calcium ferrites and the disintegration of larger grains.
The small grains also decrease in number as more reduction reactions occur.

Though the above describes just one of the possible scenarios and is quite simplified, it is consistent with the previous results that were obtained using various experimental and theoretical approaches in terms of the materials science.

Further analysis is required to fully understand the reaction process of iron ore sinters, but our analysis provides some insight into the possible scenarios and, importantly, gives us interpretable results.

Finally, we discuss the limitations of the proposed method.
This method works well when there is a one-to-one relationship between the features contained in the PDs of different dimensions.
Fortunately, our X-CT data satisfy this condition, but not all datasets do. 
We consider the following example to illustrate the problem.
Suppose the 0-dimensional PDs have two feature distributions $\psi_{0, 0}$ and $\psi_{0, 1}$, the 1-dimensional PDs also have two feature distributions $\psi_{1, 0}$ and $\psi_{1, 1}$, and the 2-dimensional PDs have one feature distribution $\psi_2$. Furthermore, suppose that there is a one-to-one correlation between $\psi_{0,0}$ and $\psi_{1,0}$, and between $\psi_{0,1}$ and $\psi_{1,1}$, but $\psi_2$ is related to both.
In this case, we find the following two concatenated feature distributions using our proposed method:
\begin{equation}
  \begin{bmatrix}
    \psi_{0,0} \\
    \psi_{1,0} \\
    \alpha\psi_2
  \end{bmatrix},
  \begin{bmatrix}
    \psi_{0,1} \\
    \psi_{1,1} \\
    \beta\psi_2
  \end{bmatrix},
\end{equation}
where $\alpha, \beta > 0$ are appropriate weights. The interpretation of $\psi_2$ in this case is more difficult than the example analysis of iron ore sinters.
In other words, if the data contain more complex correlations between features, the advantage of our method, which is good interpretability, is lost.

\acknowledgments

This work is partially supported by JSPS KAKENHI JP 19H00834, JSPS Grants-in-Aid for Transformative Research Areas (A) 20H05884 and 22H05109, and JST PRESTO JPMJPR1923.

\references


\end{document}